\documentclass[sigconf]{acmart}

\usepackage{mathtools}
\usepackage{color}

\setcopyright{rightsretained}

\begin{document}

\newcommand{\sL}{\mathcal{L}}
\newcommand{\sR}{\mathcal{R}}
\newcommand{\sH}{\mathcal{H}}
\newcommand{\sV}{\mathcal{V}}
\newcommand{\Zd}{Z_{\downarrow}}
\newcommand{\Zt}{Z_{\leftrightarrow}}
\newcommand{\rank}{\text{rank}}
\newcommand{\diag}{\text{diag}}

\title{Fast Algorithms for Displacement and Low-Rank Structured Matrices}

\author{Shivkumar Chandrasekaran}
\affiliation{%
  \institution{University of California Santa Barbara}
  \streetaddress{Department of Electrical \& Computer Engineering}
  \city{Santa Barbara}
  \state{California}
  \postcode{93106-9560}
}
\email{shiv@ucsb.edu}

\author{Nithin Govindarajan}
\affiliation{%
  \institution{University of California Santa Barbara}
  \streetaddress{Department of Mechanical Engineering}
  \city{Santa Barbara}
  \state{California}
  \postcode{93106-9560}
}
\email{ngovindarajan@ucsb.edu}

\author{Abhejit Rajagopal}
\affiliation{%
  \institution{University of California Santa Barbara}
  \streetaddress{Department of Electrical \& Computer Engineering}
  \city{Santa Barbara}
  \state{California}
  \postcode{93106-9560}
}
\email{abhejit@ucsb.edu}

\renewcommand{\shortauthors}{S.~Chandrasekaran}

\newcommand{\nithin}[1]{ {\color{red} #1} }
\newcommand{\abhe}[1]{ {\color{orange} #1} }
\newcommand{\shiv}[1]{ {\color{green} #1} }

\begin{abstract}
  This tutorial provides an introduction to the development of fast
  matrix algorithms based on the notions of displacement and various
  low-rank structures.
\end{abstract}

%

\copyrightyear{2018}
\acmYear{2018}
\setcopyright{acmlicensed}
\acmConference[ISSAC '18]{2018 ACM International Symposium on
Symbolic and Algebraic Computation}{July 16--19, 2018}{New York,
NY, USA}
\acmBooktitle{ISSAC '18: 2018 ACM International Symposium on
Symbolic and Algebraic Computation, July 16--19, 2018, New York, NY,
USA}
\acmPrice{15.00}
\acmDOI{10.1145/3208976.3209025}
\acmISBN{978-1-4503-5550-6/18/07}

\keywords{Fast matrix algorithms, displacement structure, low-rank, sequentially semi-separable,
hierarchically semi-separable, Toeplitz, Hankel, Vandermonde, Cauchy, Hilbert, Fast multipole method}

\maketitle

\section{Introduction}

In this tutorial we give a broad introduction to the class of
displacement structured and other low-rank structured matrices that
have come to recent prominence. Due to page limitations this paper only provides a short summary of the actual topics that will be explained in the tutorial itself.

As is well-known, it costs roughly
$O(n^3)$ arithmetic operations to multiply two $n\times n$ generic matrices
(as long as no Strassen-style fast multiplication algorithm is
used\footnote{We discard Strassen style algorithms as they are slower
  than the standard algorithm for reasonable values of
  $n$.}).
Similarly, the solution to a generic linear system of of $n$ equations in
$n$ unknowns requires $O(n^3)$ operations (again assuming we do not
use a fast matrix multiplication algorithm).

However, when we consider a \textit{structured} family of matrices it
might be possible to find practically fast algorithms. Famous examples
include the FFT method~\cite{cooley1965algorithm} for multiplying the discrete Fourier series
matrix rapidly. As a fallout of the FFT, we get a fast method for
convolution and as a result a large number of fast algorithms for
polynomials~\cite{GathenJoachimvonzur1999Mca}, which in turn yield fast algorithms for Toeplitz
matrices. However, many of these fast algorithms did not prove so
useful in practice~\cite{golub2013matrix}. Some of them, were too
slow for reasonable values of $n$, whereas others required large
number of bits to hold intermediate values.

Eventually some of these shortcomings were over come by exploiting
low-rank structures that are lurking in these problems (e.g.~\cite{chandrasekaran1996stabilizing,xia2012superfast,chandrasekaran2007superfast}). In this
tutorial we focus on two aspects of this. The first is the
\emph{displacement structure} approach pioneered by Kailath et. al. \cite{kailath1979displacement,kailath1995displacement}. The
second is the \emph{low-rank structured} approach pioneered by Rokhlin 
\cite{greengard1987fast,rokhlin1985rapid,rokhlin1990rapid},
Hackbusch \cite{hackbusch2015hierarchical}, Eidelman and Gohberg \cite{eidelman2014separable}, Dewilde and van der Veen \cite{DewildeVanderveen98}, and others.

\section{Displacement structure}

We eschew generality in this presentation and also do not worry about
\textit{superfast} algorithms. For simplicity we restrict ourselves to
the real field.

Let the linear operator
$\sL[A,B] : \sR^{n\times n} \rightarrow \sR^{n\times n}$ be defined as
\[
  \sL[A,B](T) = T- A T B^T,
\]
where $A, B \in \sR^{n\times n}$. We will call such a linear operator
as a displacement operator.
We say that the matrix $T$ has a displacement structure (wrt $\sL[A,B]$) if the rank of $\sL[A,B](T)$
is small compared to $n$. The standard examples are Toeplitz, Hankel,
Vandermonde, Cauchy and Pick matrices.
\begin{itemize}
\item Let $Z$ denote the down-shift matrix
\[
  Z_{i,j} = \begin{cases}
    0, & i \neq j-1 \\
    1, & i = j-1.
  \end{cases}
\]
Then it is easy to check that if $T$ is a Toeplitz matrix rank of
$\sL\left[Z,Z\right](T)$ is at most 2. Conversely, any $T$ for which
rank of $\sL[Z,Z]$ is small is called a Toeplitz-like matrix.

\item Let $x_i \in \sR$ and let $D(x) = \diag\{x_i\}$. Let $V$ be the
  Vandermonde matrix: $V_{i,j}=x_i^j$. Then rank of $\sL[D(x),Z](V)$
  is at most 1.

\item Let $0\neq y_j\in\sR$. Let $C$ be the Cauchy matrix:
  $C_{i,j}=1/(y_i-x_j)$. Then rank of $\sL[D^{-1}(y),D(x)](C)$ is at
  most 1. Conversely any matrix $C$ for which rank of $\sL[D(y), D(x)](C)$
  is small is called a generalized Pick matrix.
\end{itemize}
What is surprising is that the inverse has the same structure too. In
fact
\[
  \rank\left(\sL[A,B](T)\right) =
  \rank\left(\sL[B,A]\left(T^{-1}\right)\right).
\]
More trivially
\begin{align*}
  \rank\left(\sL[A,B](T_1+T_2)\right) \leq
  \rank\left(\sL[A,B](T_1)\right) + \\ \rank\left(\sL[A,B](T_2)\right).
\end{align*}
Unfortunately the result for products is not as nice:
\begin{align*}
  \rank\left(\sL[A,B](T_1T_2)\right) \leq
  \rank\left(\sL[A,B](T_1)\right) + \\
  \rank\left(\sL[A,B](T_2)\right) +
  \rank\left(B^TA-I\right).
\end{align*}
In particular this is rather disastrous for generalized Pick (Cauchy)
matrices, which will be rectified later.

We next look for fast algorithms that can exploit the displacement
structure. Suppose that
\[
  T - A T B^T= P Q^T,
\]
where $P, Q \in \sR^{n \times p}$. If $\sL[A,B]$ is invertible then
the pair $(P,Q)$ can be used as a more efficient representation of the
matrix $T$ as we would require only $O(n p)$ numbers versus the usual
$O(n^2)$ numbers. In this case the pair $(P,Q)$ will be called
\textbf{generators} for $T$. The question is: can we carry out the
usual matrix operations faster using the pair $(P,Q)$? Surprisingly,
this question is easier to answer for $T^{-1}x$ rather than for $Tx$,
for $x\in\sR^n$. In particular we will show that the $L U$
factorization (Gaussian elimination) of $T$ can be computed quickly in
$O(n^2)$ arithmetic operations rather than the standard $O(n^3)$
operations, provided inverting $\sL[A,B]$ itself is cheap.

We note that the displacement rank is invariant under similarity
transformations
\[
  \rank\left(\sL[VAV^{-1},WBW^{-1}](VTW^T)\right) =
  \rank\left(\sL[A,B](T)\right).
\]
Therefore it is convenient to assume that $A$ and $B$ are lower
triangular matrices from now on (say using the Schur decomposition). Note that for this to be effective in
fast algorithms, the computation and application of $V$, $V^{-1}$ $W$
and $W^{-1}$, must cost less than $O(n^3)$ operations. There are
extensions of the method to the case when $A$ or $B$ is a lower
Hessenberg matrix \cite{heinig2001schur}, but we do not cover it here.

Let $\sL[A,B](T)$ have low-rank with $A$ and $B$ lower triangular and
let
\[
  T = \bordermatrix{ & 1 & n-1 \cr
    1 & T_{11} & T_{12} \cr
    n-1 & T_{21} & T_{22} }, \qquad S = T_{22}- T_{21}T_{11}^{-1}T_{12},
\]
where $S$ is the Schur complement of $T$ and we have have assumed that
$T_{11} \neq 0$. The key idea of fast Gaussian elimination is that $S$
has a low-rank displacement structure that can be quickly computed
from that of $T$. To see this, let
\[
  A = \bordermatrix{ & 1 & n-1 \cr
    1 & A_{11} & 0 \cr
    n-1 & A_{21} & A_{22} }, \qquad
  B = \bordermatrix{ & 1 & n-1 \cr
    1 & B_{11} & 0 \cr
    n-1 & B_{21} & B_{22} }.
\]
Then it can be readily shown that rank of $\sL[A_{22},B_{22}](S)$ is
at most rank of $\sL[A,B](T)$, and furthermore the generators for $S$
can be computed in $O(np^2)$ operations provided that the first column
and row of $T$ can be generated in $O(n)$ operations\footnote{This
  requires that $A$ and $B$ be specially structured lower-triangular
  matrices. The shift-down and diagonal matrices are good examples.}
from its generators $(P,Q)$. The recursive application of this idea
leads to the fast generalized Schur algorithm for computing the $LU$
factorization of $T$ in $O(n^2)$ operations.

For many applications already this speed-up is sufficient. However
there is sometimes a need for more. For example if $T$ is a Toeplitz
matrix we can compute $Tx$ in $O(n \log_2 n)$ operations using the
FFT. The Gohberg-Semencul formulas \cite{GohbergSemencul72,chun1989constructive} show that there is a representation
for $T^{-1}$ (which might require some time to compute) such that
$T^{-1}x$ can also be computed in $O(n \log_2 n)$ operations via the
FFT. The displacement structure approach can also shed light on this
situation.

First note that $T^{-1}$ is a Schur complement of
\[
  M = \begin{bmatrix}
    T & I \\ -I & 0
  \end{bmatrix}.
\]
So we could compute the generators of $T^{-1}$ by running the
generalized Schur algorithm
on $M$ half-way through. This would be
fast if $M$ had a short displacement rank for an appropriate
$\sL[A_M,B_M]$. One potential choice is
\[
  A_M = \begin{bmatrix}
    A & \\ & A
  \end{bmatrix}, \qquad B_M = \begin{bmatrix}
    B & \\ & B
  \end{bmatrix},
\]
provided rank of $I - A B^T$ is small. Another possibility is
\[
  A_M = \begin{bmatrix}
    A & \\ & B^{T\dagger}
  \end{bmatrix}, \qquad B_M = \begin{bmatrix}
    B & \\ & A^{T\dagger}
  \end{bmatrix},
\]
provided the orthogonal projectors $I-BB^{\dagger}$ and $I-AA^\dagger$
have small rank.

For generalized Pick matrices $P$ where the rank of $\sL[D_1,D_2](P)$
is small, with $D_i$ a diagonal matrix, we have a problem if
$\sL[D_1,D_2]$ is not invertible. In this case we can still make
progress via a rank-1 perturbation of the form $\sL[D_1+uu^T,D_2](P)$
for a well-chosen column $u$, and exploiting the fact the eigenvector
matrix of $D_1+uu^T$ is an orthogonal Cauchy matrix, which can be
multiplied rapidly using one of several techniques. This rank 1
perturbation will increase the displacement rank of $P$ by at most 1 so
the fast generalized Schur algorithm can still be deployed.

Once the generators for $T^{-1}$ have been computed, the question
naturally arises as to how to multiply $T^{-1}x$ quickly. If
$|\lambda(A)| < 1$ and $|\lambda(B)| < 1$, then we have the series
solution for $\sL[A,B](T) = PQ^T$ as:
\begin{equation}
  T = \sum_{l=0}^\infty \left(A^lP\right) \left(B^lQ\right)^T.
  \label{eqn:krylov}
\end{equation}
If $T$ is the inverse of a Toeplitz matrix and $A=B=Z$, then this
gives the celebrated Gohberg--Semencul formula as $A$ and $B$ are
nilpotent. When $A$ and $B$ are diagonal then $T$ is a generalized
Cauchy (Pick) matrix and there are fast algorithms for matrix--vector
multiply~\cite{GohbergOlshevsky94}. Some generalizations are available
in~\cite{BostanJMS17}. However there seems to be no general approach
for multiplying an arbitrary low displacement-rank matrix with a
vector for a general class of \textit{nice} $A$ and $B$.

We do make the general observation that if the sum~(\ref{eqn:krylov})
is converging rapidly then one can easily construct a fast
\textit{approximate} matrix--vector multiplication algorithm. But
those matters are better dealt with later in this tutorial.

It is easy to see how to compute the generators for $T_1+T_2$ quickly
from those of $T_i$ provided they have low displacement-rank for the
same displacement operators $\sL[A,B]$. Similarly it is not difficult
to quickly compute the generators of $T_1T_2$ provided that we can
rapidly multiply by both $A$ and $B$, and that a low rank
factorization of $I-B^TA$ is quickly available.

The big draw back of displacement structure approach is that there has
been no progress on generalizing it to nested structures like
Toeplitz-block-Toeplitz matrices. This might in turn be related to the
fact that the inverses of higher-order Sylvester--Stein operators are
not well-understood.

\section{Sequentially Semi-Separable (SSS) representations}

Just like displacement structure theory was born from \textit{systems}
theory (study of linear time-invariant [LTI] systems), another study
of low-rank structured matrices was born from trying to generalize
systems theory to time-varying systems~\cite{DewildeVanderveen98}. We
do not attempt a generalized approach in this tutorial. (More information can be found in~\cite{chandrasekaran2005some,chandrasekaran2007fast,chandrasekaran2007superfast}.)

For any square matrix $A$ consider partitions of the form
\[
  A = \bordermatrix{ & m & n-m \cr
    m & A_{11} & A_{12} \cr
    n-m & A_{21} & A_{22} }.
\]
We will call such off-diagonal blocks as $A_{12}$ and $A_{21}$ as
Hankel blocks. Note that these off-diagonal blocks never cross the
principal diagonal and always extend as far to the corners as
possible. In this section we study families of matrices for which the
ranks of all Hankel blocks are small compared to the matrix size
$n$. In particular it turns out that there is a non-linear
representation of the matrix that captures this structure precisely
and permits a full spectrum of fast linear (in matrix size $n$)
algorithms.

The key question is what constraints are placed on the matrix entries
by the requirement that two overlapping Hankel blocks must have low
rank? Consider the following partition of $A$:
\[
  A = \begin{bmatrix}
    A_{11} & A_{12} & A_{13} \\
    A_{21} & A_{22} & A_{23} \\
    A_{31} & A_{32} & A_{33} \\
  \end{bmatrix}
\]
where both $\sH_1 = \begin{bmatrix}A_{12} & A_{13} \end{bmatrix}$ and
$\sH_2 = \begin{bmatrix} A_{13}^T & A_{23}^T \end{bmatrix}^T$ are
Hankel blocks. From the presence of the shared common block $A_{13}$
we see that the column space of $\sH_1$ must be related to the column
space of $\sH_2$, and similarly the row-space of $\sH_2$ must be
related to the row-space of $\sH_1$. To that end, let:
\[
  \sH_1 = U_1 \begin{bmatrix} V_{11}^T & V_{12}^T
  \end{bmatrix}, \qquad \sH_2 = \begin{bmatrix}
    U_{21} \\ U_{22}
  \end{bmatrix} V_2^T,
\]
be full column-rank factorizations of $\sH_i$ conformally partitioned
with $A_{13}=U_1V_{12}^T=U_{21}V_2^T$. For a fixed choice of the two
factorizations there is a unique matrix $W$, which we call a
transition operator, such that
\[
  U_{21} = U_1 W, \qquad V_{12}^T=W V_2^T.
\]
If the rank of $\sH_i$ is $r_i$ then $W \in \sR^{r_1\times r_2}$.

If we now choose $\sum_{i=1}^p n_i = n$, and use it to partition $A$
into a block $p\times p$ matrix, such that the sub-block
$A_{ij} \in \sR^{n_i\times n_j}$, then we can use the above idea
repeatedly on adjacent Hankel blocks to construct a representation of
the blocks of $A$ of the following form:
\[
  A_{ij} = \begin{cases}
    D_i, & i = j, \\
    U_i W_{i+1} W_{i+2} \cdots W_{j-1} V_j^T, & i < j, \\
    P_i R_{i-1} R_{i-2} \cdots R_{j+1} Q_j^T, & i > j.
  \end{cases}    
\]

A simple example is a banded matrix which clearly has the desired
low-rank property for the Hankel blocks. In this case if we choose
$n_i$ to be equal to the band-width we can write down the components
directly.
\begin{eqnarray*}
  D_i & = & A_{i,i} \\
  U_i &=& A_{i,i+1} \\
  V_i & = & I \\
  P_i &=& A_{i,i-1} \\
  Q_i &=& I \\
  R_i &=& 0 \\
  W_i &=& 0. \\
\end{eqnarray*}
To understand this representation better we look at how you can
compute $Ax=b$ fast. A little bit of algebra reveals that the
following recursions will do the job:
\begin{eqnarray}
  g_i & = & V_i^Tx_i + W_i g_{i+1} \nonumber \\
  h_i & = & Q_i^Tx_i + R_i h_{i-1}  \nonumber \\
  b_i & = & D_i x_i + U_i g_{i+1} + P_i h_{i-1}, \label{eqn:sssrec}
\end{eqnarray}
where all undefined variables are assumed to be empty matrices of
suitable size, and $x_i$ and $b_i$ are conformal partitions of $x$ and
$b$ respectively. It is clear that if all Hankel block ranks are small
compared with the matrix size then this is a linear time algorithm for
computing $Ax$.

However there are more to these recursions. We make the notation, for
example, that
\[
  D = \begin{bmatrix}
    D_1 & \\
    & D_2 \\
    & & \ddots \\
    & & & D_p
  \end{bmatrix} \qquad g = \begin{bmatrix}
    g_1 \\ g_2 \\ \vdots \\ g_p
  \end{bmatrix}.
\]
The recursions can then be written in matrix form as
\begin{equation}
  \begin{bmatrix}
    I-W Z^T & 0 & -V^T \\
    0 & I - R Z & -Q^T \\
    UZ^T & PZ & D
  \end{bmatrix} \begin{bmatrix}
    g \\ h \\ x
  \end{bmatrix} = \begin{bmatrix}
    0 \\ 0 \\ b
  \end{bmatrix}. \label{eqn:sparseSSS}
\end{equation}
From this we get the \textit{diagonal} representation of $A$ as
\begin{equation}
  A = D + U Z^T(I-WZ^T)^{-1}V^T + P Z(I - RZ)^{-1}Q^T. \label{eqn:sssdiag}
\end{equation}
In other words $A$ is just the Schur complement of a large sparse
matrix. In particular we will call $(D, U, W, V, P, R, Q)$ as the SSS
representation of $A$.

Furthermore there is a re-ordering of the unknowns in
equation~(\ref{eqn:sparseSSS}) into the sequence $(g_i, h_i, x_i)$,
and if the rows of the sparse matrix in~(\ref{eqn:sparseSSS}) are also
re-ordered the same way we can then see that the the sparse matrix
will be a linear graph by referring to the original recursions
in~(\ref{eqn:sssrec}). Therefore it follows that in this ordering the
sparse matrix has a no fill-in elimination order and that therefore
there is a fast algorithm (sparse Gaussian elimination) to quickly
compute $x$ (and incidentally also $g$ and $h$) from $b$. Some
elementary calculations then reveal that the resulting algorithm wll
be linear in the matrix size $n$.

However more operations involving SSS representations can be done in
linear time. The key to understanding this is the simple fact that
low-rank Hankel blocks imply short SSS representations and vice versa.

For example, since
\[
  \begin{bmatrix}
    A_{11} & A_{12} \\
    A_{21} & A_{22}
  \end{bmatrix}   \begin{bmatrix}
    B_{11} & B_{12} \\
    B_{21} & B_{22}
  \end{bmatrix} =
  \begin{bmatrix}
    * & A_{11}B_{12}+A_{12}B_{22} \\
    * & *
  \end{bmatrix}
\]
it follows that the product of two SSS matrices will be another SSS
matrix whose Hankel block ranks are at most the sum of the
corresponding Hankel blocks of the original SSS matrices. With some
effort a linear time algorithm can be constructed for producing
the SSS representation of $AB$ given those of $A$ and $B$.

Similarly from
\[
    \begin{bmatrix}
    A_{11} & A_{12} \\
    A_{21} & A_{22}
  \end{bmatrix}^{-1} =
    \begin{bmatrix}
    * & *A_{12}* \\
    *A_{21}* & *
  \end{bmatrix},
\]
where the $*$'s denote matrices that we do not care about, it follows
that the Hankel block ranks of the inverse are at most those of the
original matrix. It therefore follows, for example, that both banded
matrices and their inverse have short SSS representations. Though not
all matrices with short SSS representations are the inverses of banded
matrices, the diagonal representation~(\ref{eqn:sssdiag}) shows that
every SSS representation \textit{is} the Schur complement of a larger
banded matrix in the right ordering.

One way to quickly compute the SSS representation of the inverse is as
follows. First we observe that since
\[
    \begin{bmatrix}
    A_{11} & A_{12} \\
    A_{21} & A_{22}
  \end{bmatrix} =    \begin{bmatrix}
    I & 0 \\
    A_{21}A_{11}^{-1}& I
  \end{bmatrix}  \begin{bmatrix}
    A_{11} & A_{12} \\
    0 & A_{22}-A_{21}A_{11}^{-1}A_{12}
  \end{bmatrix},
\]
both $L$ and $U$ in the $LU$ factorization of $A$ will have exactly
the same Hankel block ranks as $A$. From this with some algebra one
can construct a linear time algorithm to compute the SSS
representations of $L$ and $U$ from that of $A$. Then, using the
diagonal representation and the Woodbury formula, we can compute in linear time the SSS
representation of $L^{-1}$ and $U^{-1}$. We can then use the aforementioned linear time multiplication algorithm to find the SSS
representation of $U^{-1}L^{-1} = A^{-1}$.

We also note that there are linear time algorithms for the $ULV$
factorization of $A$ in SSS form, if growth factor becomes an issue.

Clearly the sum of two SSS matrices is another SSS matrix. So
effectively matrix algebra is fast in the SSS representation and this
has proven to be tremendously useful in practice.

\section{Hierarchically Semi-Separable (HSS) representations}
There is another non-linear matrix representation that is also capable
of providing a linear time matrix algebra. This was born independently
from efforts to speed up the application of integral operators that
arise in potential theory in the fundamental work of Greengard and
Rokhlin on the fast multi-pole method (FMM). The representation we
present now (HSS) is a special case of this more general class of FMM
representations. The more general class is not \textit{closed} under
inversion and multiplication, but the HSS representation is. More detailed information can be found in~\cite{chandrasekaran2006fast,lessel2016fast,chandrasekaran2006fast2,dewilde2006hierarchical,sheng2007algorithms,lyons2005fast,Timpals,starr}.

The HSS representation uses a slightly different definition of Hankel
blocks. Let
\[
  A = \begin{bmatrix}
    A_{11} & A_{12} & A_{13} \\
    A_{21} & A_{22} & A_{23} \\
    A_{31} & A_{32} & A_{33} \\
  \end{bmatrix}.
\]
Then we call $\begin{bmatrix} A_{21} & A_{23}
\end{bmatrix}$ a row Hankel block, and we call $\begin{bmatrix}
  A_{12}^T & A_{32}^T
\end{bmatrix}^T$ a column Hankel block. The HSS representation
exploits the low-rank structure of these types of blocks, but not
\textit{all} of them.
It restricts itself to a fixed set hierarchical set of blocks instead. (This is to be contrasted with the SSS representation which does capture the low-rank structure of all the relevant Hankel blocks.)

To this end we need some notation. Let $A=A_{0;0,0} \in \sR^{n\times n}$, and let $n=n_{0;0}$. We assume that there is a partition tree associated to $A$, which is defined as follows. We assume that for each $0\leq k< K$ there are non-negative integers $n_{k;i}$, for $i=0,\ldots,2^k-1$, such that $n_{k;i}=n_{k+1;2i}+n_{k+1;2i+1}$. Note that these numbers can be naturally associated with a binary tree which we call the partition tree. (In general we do not need to have a complete binary tree, but we ignore that generalization here.)

Based on this partition tree we recursively partition $A$ as follows:
\[
    A_{k;i,j} = \bordermatrix{ & n_{k+1;2j} & n_{k+1;2j+1} \cr
    n_{k+1;2i} & A_{k+1;2i,2j} & A_{k+1;2i,2j+1} \cr
    n_{k+1;2i+1} & A_{k+1;2i+1,2j} & A_{k+1;2i+1,2j+1}},
\]
and these sub-blocks can be viewed as edges on the partition tree.

We define the row Hankel blocks at each level as:
\[
\sH_{k;i} = \begin{bmatrix}
A_{k;i,0} & \cdots & A_{k;i,i-1} & A_{k,i,i+1} & \cdots & A_{k,i,2^k-1}
\end{bmatrix}.
\]
That is, $\sH_{k;i}$ is the i'th row at the k'th level with the block $A_{k;i,i}$ deleted. In analogous manner, one can define column Hankel blocks $\sV_{k;i}$, where the block columns of $A$ at the $k$-th level are considered instead. Next let the matrix $V_{k;i}$ be such that its columns form a basis for the row space of $\sH_{k;i}$. Similarly let the matrix $U_{k;i}$ be such that its columns form a basis for the column space of $\sH_{k;i}$.

Since $\sH_{k,i}$ shares sub-matrices with $\sH_{k+1;2i}$ and $\sH_{k+1;2i+1}$, it follows that there are \textit{translation} matrices $R_{k;i}$ such that
\[
    U_{k;i} = \begin{bmatrix}
    U_{k+1;2i} R_{k+1;2i} \\
    U_{k+1;2i+1} R_{k+1;2i+1}
    \end{bmatrix}.
\]
Similarly there are translation matrices $W_{k;i}$ associated with $V_{k;i}$. 
We define the \textit{expansion coefficients} $B_{k;i,j}$ as follows:
\begin{eqnarray*}
    A_{k;2i,2i+1} &=& U_{k;2i} B_{k;2i,2i+1} V_{k;2i+1}^T \\
    A_{k;2i+1,2i} &=& U_{k;2i+1} B_{k;2i+1,2i} V_{k;2i}^T,
\end{eqnarray*}
for $0<k\leq K$. We also define
\[
    A_{K;i,i} = D_{K;i}.
\]

We next observe that the translation matrices can be much smaller than the basis matrices. Therefore in the HSS representation we store only $U_{K;i}$ and $V_{K;i}$. To recover the other basis matrices we use $R_{k;i}$ and $W_{:k;i}$ instead. Therefore the complete HSS representation of $A$ is the partition tree $n_{k;i}$ along with the leaf-level matrices $D_{K;i}$, $U_{K;i}$, $V_{K;i}$, the translation matrices $R_{k;i}$, $W_{k;i}$, and the expansion coefficients $B_{k;2i,2i+1}$, $B_{k;2i+1,2i}$. With this we can check that every entry of $A$ can be uniquely recovered from the HSS representation.

If $p_{k;i}$ denotes the rank of $\sH_{k;i}$ and $q_{k;i}$ the rank of $\sV_{k;i}$, then it can be verified, for example, that $R_{k+1;2i} \in \sR^{p_{k+1;2i}\times p_{k;i}}$ and $B_{k;2i,2i+1}\in \sR^{p_{k;2i}\times q_{k;2i+1}}$. So the ranks of the row and column Hankel blocks determines how small the HSS representation is.

Just as in the SSS case, there exists an $O(n^2)$ algorithm to construct an optimal HSS representation directly from the entries of the matrix~\cite{chandrasekaran2006fast2}. However, in special cases we can do much better. For sparse matrices the construction can be done in linear time~\cite{chandrasekaran2006fast2}. For matrices of the form $A_{i,j}=f(x_i,y_j)$, where $f:\sR^d\rightarrow\sR$ the FMM techniques of Greengard and Rokhlin can be used to compute the HSS representation of $A$ in linear time provided $f$ satisfies some nice properties~\cite{greengard1987fast}. However, in practice the HSS representation is most often computed quickly from the fact that matrix algebra in HSS form can be done quickly.

Just as in the SSS case, the key observation is that there is a short HSS representation as long as the row and column Hankel blocks have small rank. Then one observes for example, that the sum of two HSS matrices will have column (row) Hankel block ranks that are the sum of the corresponding column (row) Hankel blocks of the summands. Therefore it is not surprising that there is a linear time algorithm to add two HSS matrices.

The superfast multiplication of two HSS matrices is based on the following observation:
\begin{align*}
    \begin{bmatrix}
    A_{11} & A_{12} & A_{13} \\
    A_{21} & A_{22} & A_{23} \\
    A_{31} & A_{32} & A_{33} \\
  \end{bmatrix} \begin{bmatrix}
    B_{11} & B_{12} & B_{13} \\
    B_{21} & B_{22} & B_{23} \\
    B_{31} & B_{32} & B_{33} \\
  \end{bmatrix} = \\
  \begin{bmatrix}
    * & *B_{12}+A_{12}*+*B_{32} & * \\
    A_{21}*+*B_{21}+A_{23}* & * & A_{21}*+*B_{23}+A_{23}* \\
    * & *B_{12}+A_{32}*+*B_{32} & * \\
  \end{bmatrix},
\end{align*}
where $*$'s denote matrices whose ranks are irrelevant.  From this we see that the ranks of the column (row) Hankel blocks of the product are the sums of the ranks of the corresponding colum (row) Hankel blocks of the multiplicands. Again there is a linear time algorithm to compute the HSS representation of the product from the HSS representation of the multiplicands.

To give a flavor of how HSS algorithms are constructed we go back to the problem of fast matrix-vector multiplications $Ax=b$, where $A$ is in HSS form. We will use the notation $x_{K;i}$ to denote the $i$-block of $x$ when it is partitioned according to the cuts at the leaf level of the partition tree.
Then with a little bit of algebra one can show that the following recursions will do the job:
\begin{eqnarray*}
  g_{K;i} &=& V_{K;i}^T x_{K;i} \\
  g_{k;i} &=& W_{k+1;2i}^Tg_{k+1:2i} + W_{k+1;2i+1}^Tg_{k+1;2i+1}, \qquad k < K, \\
  f_{0;0} &=& [] \\
  f_{k;2i} &=& R_{k;2i} f_{k-1;i} + B_{k;2i,2i+1}g_{k;2i+1}, \qquad 0<k\leq K, \\
  f_{k;2i+1} &=& R_{k;2i+1} f_{k-1;i} + B_{k;2i+1,2i}g_{k;2i}, \qquad 0<k\leq K, \\
  b_{K;i} &=& D_{K;i} x_{K;i} + U_{K;i} f_{K;i}.
\end{eqnarray*}
It is a clear that this leads to a linear time algorithm. Just like in the SSS case these recursions can also be used to give a diagonal representation of $A$ and to also give a linear time solver for computing $x$ given $b$ via sparse Gaussian elimination.

Towards this define the pair, $\Zd$ and $\Zt$, of linear operators on the binary partition tree, via the equations:
\begin{eqnarray*}
  \left(\Zd x\right)_{k;i} &=& x_{k-1; \left\lfloor \frac{i}{2} \right\rfloor} \\
  \left(\Zt x\right)_{k;2i} &=& x_{k;2i+1} \\
  \left(\Zt x\right)_{k;2i+1} &=& x_{k;2i}.
\end{eqnarray*}
Also define the linear projection operator $P$ on the binary partition tree such that $P^Tg$ restricts $g$ just to the leaves of the tree, $(P^Tg)_i=g_{K;i}$.

As before also define the block diagonal matrices $D$, $U$ and $V$, which only consist of the entries $D_{K;i}$, $U_{K;i}$ and $V_{K;i}$ that are on the leaves of the binary partition tree. Also define the block diagonal matrices $R$ and $W$ that consist of the entries $R_{k;i}$ and $W_{k;i}$ which are defined on all nodes of the binary partition tree. Finally define the block diagonal matrix $B$ which consists of the entries $B_{k;2i,2i+1}$ and $B_{k;2i+1,2i}$ by making a natural association with the nodes of the binary tree. Similarly we will use $g$ to denote a column vector containing all the $g_{k;i}$'s.

Then we can re-write the fast recursions for the multiplication $Ax$ as
\begin{eqnarray*}
g &=& \Zd^T W^T g + P V^T x \\
f &=& R \Zd f + B \Zt g \\
b &=& D x + U P^T f.
\end{eqnarray*}
In matrix form this appears as:
\[
\begin{bmatrix}
    I - \Zd^T W^T & 0 & -P V^T \\
    -B \Zt & I - R \Zd & 0 \\
    0 & U P^T & D
\end{bmatrix} \begin{bmatrix}
    g \\ f \\ x
\end{bmatrix} = \begin{bmatrix}
    0 \\ 0 \\ b
\end{bmatrix}.
\]
This implies that the HSS diagonal representation of a matrix $A$ is given by
\[
    A = D + U P^T (I-R \Zd)^{-1} B \Zt (I - \Zd^T W^T)^{-1} P V^T,
\]
which also shows that $A$ is just the Schur complement of a larger block sparse matrix with the binary partition tree with edges \textit{only} between siblings, as the incidence graph.

Since such a binary partition tree has a fill-in free elimination order, this also shows that we can get a linear time solver for constructing $x$ given $b$, by first re-ordering the unknowns (and the equations) in the order: $(g_{K;i}, f_{K;i}, x_{K;i})$, followed by $(g_{k;i}, f_{k;i})$.

If numerical stability is a concern, it is also clear that a linear time sparse $QR$ factorization algorithm can also be constructed by similar considerations.

If $A$ has a short HSS representation then so does $A^{-1}$ and this can be computed in linear time from that of $A$. To see how and why this is possible, first observe that:
\begin{align*}
  \begin{bmatrix}
    A_{11} & A_{12} & A_{13} \\
    A_{21} & A_{22} & A_{23} \\
    A_{31} & A_{32} & A_{33} \\
  \end{bmatrix} = \\
  \begin{bmatrix}
      I & 0 & 0 \\
      A_{21} * & I & 0 \\
      * & (A_{32}+*A_{12})* & I
  \end{bmatrix} \begin{bmatrix}
    * & A_{12} & * \\
    0 & * & A_{23}+A_{21}* \\
    0 & 0 & * \\
  \end{bmatrix},
\end{align*}
where again $*$'s denote matrices whose ranks are irrelevant.
From this it follows that the $LU$ factors of $A$ will have the same column and row Hankel block ranks as that of $A$ itself. Furthermore, with a little bit bit of algebra, a set of fast linear time recursions can be worked out for computing the HSS representations of the $LU$ factors from that of $A$.

From:
\begin{align*}
  \begin{bmatrix}
    U_{11} & U_{12} & U_{13} \\
    0 & U_{22} & U_{23} \\
    0 & 0 & U_{33} \\
  \end{bmatrix}^{-1} =
  \begin{bmatrix}
      * & *U_{12}* & * \\
      0 & * & *U_{23}* \\
      0 & 0 & *
  \end{bmatrix},
\end{align*}
we observe that the columns and row Hankel blocks of the inverse of a triangular matrix have the same ranks as those of the original matrix. Furthermore a linear time algorithm can be devised for computing the HSS representation of the inverse of a triangular HSS matrix.

We already know that the product of two HSS matrices will have short HSS representations, but the Hankel block ranks will add up. So when we multiply $U^{-1}L^{-1}=A^{-1}$ there is the danger that we will end up with larger column and row Hankel block ranks. However, we note that:
\begin{align*}
\begin{bmatrix}
    A_{11} & A_{12} & A_{13} \\
    A_{21} & A_{22} & A_{23} \\
    A_{31} & A_{32} & A_{33} \\
  \end{bmatrix}^{-1} = \\
  \begin{bmatrix}
      * & (*A_{12}+*A_{32})S_v & * \\
      S_h(A_{21}*+A_{23}*) & * & S_h(A_{21}*+A_{23}*) \\
      * & (*A_{12}+*A_{32})S_v & *
  \end{bmatrix},
\end{align*}
which shows that the column and row Hankel block ranks of $A^{-1}$ are at most that of $A$ itself. Therefore with a little algebra we can devise linear time algorithms to compute the HSS representation of $A^{-1}$ from those of the $LU$ factors of $A$. The above formula is also the reason for the definition of row and column Hankel blocks.

\section{Conclusion}

We have presented a basic outline of the displacement structured matrices and the construction of fast solvers for them. However the biggest open question in this area is whether this approach generalizes to dealing with Toeplitz--block--Toeplitz (TBT) matrices?

Currently the most efficient solvers convert Toeplitz matrices to Cauchy-like matrices via the FFT and exploit their HSS structure instead. However, even this approach does not seem to extend to TBT matrices, though in general the HSS approach does not do too badly if there is an underlying 2D kernel function (via the FMM representation).

For SSS and HSS type representations the biggest open question is the exact rank structure of the inverse of discrete 2D Laplace like matrices. There are several papers on the approximate low-rank structure in special cases, but the general question remains quite open.

We have made some remarks recently on this problem~\cite{chandrasekaran2018HK}. Our key observation is to note that the complexity of SSS algorithms depends critically on the underlying linear graph, while that of HSS depends on the special binary partition tree. Conversely the lack of a fast direct solver for FMM representations can be traced to the difficulty of doing Gaussian elimination quickly on a the more complicated FMM tree. This immediately raises the issue of tying the graph structure more intimately to the matrix representation and the associated fast algorithms.

This problem is important as our understanding of the inverse of discrete 3D Laplace like matrices is not sharp enough to yield fast enough  practical solvers and the underlying graph is the 3D mesh.


\bibliographystyle{plain}
\bibliography{sample-bibliography}

\end{document}